\documentclass[11pt]{article}
\usepackage{amsfonts}
\usepackage{graphics}
\usepackage{indentfirst}
\usepackage{cite}
\usepackage{latexsym}
\usepackage{amsmath}
\usepackage{amssymb}
\usepackage[dvips]{epsfig}
\usepackage{amscd}

\newtheorem{theorem}{Theorem}[section]
\newtheorem{remark}{Remark}[section]

\newtheorem{definition}{Definition}[section]
\newtheorem{lemma}[theorem]{Lemma}

\newtheorem{proposition}[theorem]{Proposition}

\newcommand{\n}{\rho}
\newcommand{\nr}{\tilde\rho^R}
\newcommand{\ti}{\tilde}
\newcommand{\mr}{\mathbb{R}}

\def\pf{{\it Proof.}  }
\renewcommand{\div}{ {\rm div }  }

\newcommand{\na}{\nabla }
\newcommand{\vp}{\varphi }
\newcommand{\pa}{\partial}

\newcommand{\bt}{\begin{theorem}}
\newcommand{\bl}{\begin{lemma}}
\newcommand{\el}{\end{lemma}}
\newcommand{\et}{\end{theorem}}

\newcommand{\de}{\delta}
\newcommand{\ve}{\varepsilon}
\newcommand{\la}{\label}

\newcommand{\bn}{\begin{eqnarray}}
\newcommand{\en}{\end{eqnarray}}
\newcommand{\bnn}{\begin{eqnarray*}}
\newcommand{\enn}{\end{eqnarray*}}

\newcommand{\bnnn}{\begin{eqnarray*}}
\newcommand{\ennn}{\end{eqnarray*}}
\newcommand{\ben}{\begin{enumerate}}
\newcommand{\een}{\end{enumerate}}

\newcommand{\xl}{\left}
\newcommand{\xr}{\right}

\newcommand{\ba}{\begin{aligned}}
\newcommand{\ea}{\end{aligned}}
\newcommand{\be}{\begin{equation}}
\newcommand{\ee}{\end{equation}}

\def\p{\partial}
\def\norm[#1]#2{\|#2\|_{#1}}

\def\lap{\triangle}

\def\ep{\varepsilon}

\def\rr{\mathbb{R}^2}
\def\O{B_R}

\makeatletter
\@addtoreset{equation}{section}
\makeatother

\title{On Local Strong Solutions to the Cauchy Problem of Two-Dimensional
Density-Dependent Magnetohydrodynamic Equations with Vacuum
\thanks{B. L\"u is  supported by  NNSFC Tianyuan (No. 11426131) and Natural Science Foundation of Jiangxi Province (No. 20142BAB211006).
}
}
\date{}
\author{ Boqiang L\"u\thanks{College of Mathematics and Information Science, Nanchang Hangkong University, Nanchang 330063, P. R.  China({\tt lvbq86@163.com}). }
\quad   Zhonghai Xu\thanks{College of Science, Northeast Dianli University, Jilin 132013, P. R. China ({\tt xuzhonghai@163.com}).
 }
\quad Xin Zhong\thanks{Corresponding author. Institute of Applied Mathematics, AMSS,
Chinese Academy of Sciences, Beijing 100190, P. R. China({\tt xzhong1014@amss.ac.cn}).
 }
 }

\begin{document}
\maketitle

\begin{abstract}  This paper concerns the   Cauchy problem of the nonhomogeneous incompressible magnetohydrodynamic (MHD) equations on the whole two-dimensional (2D) space with  vacuum as far field density. In particular, the initial density can have compact support.
 We prove that  the 2D  Cauchy problem  of the nonhomogeneous incompressible  MHD equations admits a unique local strong solution
  provided the initial density and the initial  magnetic decay not too slow at infinity.
\end{abstract}


\textbf{Math Subject Classification}: 35Q35; 76D03; 76W05.

\section{Introduction and main results}
We consider the two-dimensional nonhomogeneous incompressible  magnetohydrodynamic equations which read as follows:
\be\la{mhd}
\begin{cases}\n_t + \div(\n u) = 0,\\
 (\n u)_t + \div(\n u\otimes u) + \nabla P = \mu\Delta u + H\cdot\na H-\frac{1}{2}\na |H|^2,\\
 H_t-\nu\Delta H+u\cdot\na H-H\cdot\na u=0,\\
 \div u= \div H=0,
\end{cases}
\ee
where $t\ge 0 $ is time, $ x=(x_1,x_2)\in \Omega\subset \rr$ is the spatial coordinate, and $\rho=\rho(x,t)$, $ u =(u^1,u^2)(x,t)$, $ H =(H^1,H^2)(x,t)$, and $P=P(x,t)$ denote the density, velocity, magnetic, and pressure of the fluid, respectively; $\mu>0$ stands for the viscosity constant.
The constant $\nu>0$ is the resistivity coefficient which is inversely proportional to the electrical conductivity constant and acts as the magnetic diffusivity of magnetic fields.

Let $\Omega=\rr$ and we consider the Cauchy problem for \eqref{mhd} with $(\rho, u, H)$ vanishing at infinity (in some weak sense) and the initial conditions:
\be \la{n4}
\n(x,0)=\n_0(x), \quad  \n u (x,0)= \n_0u_0(x),\quad  H(x,0)= H_0(x),\quad x\in \Omega,
\ee
 for given initial data $\rho_0, u_0$ and $H_0$.

Magnetohydrodynamics studies the dynamics of electrically conducting fluids and the theory of the macroscopic interaction of electrically conducting fluids with a magnetic field.
In particular, if there is no electromagnetic effect,
that is $H=0$, the MHD system reduces to the Navier-Stokes equations,
which have been discussed by many mathematicians, please see  \cite{AK1973,AKM1990,CK2003,HW2014,HW2015,K1974,L1996,L2015,
F2004,H19951,L1998,HL20132,HLX2012,lx,LL2014,LX2014}  and references therein. Since the fluid motion and the magnetic field are strongly couplied and interplay interaction with each other, it is rather complicated to investigate the MHD system. Now, we briefly recall some results
concerning with the multi-dimensional nonhomogeneous incompressible MHD equations which are more relatively to our problem. Gerbeau-Le Bris \cite{leb1} and Desjardins-Le Bris \cite{leb2} studied the global existence of weak solutions with finite energy on 3D bounded domains  and on the torus, respectively. In the absence of vacuum, Abidi-Hmidi \cite{abidi1} and Abidi-Paicu \cite{abidi2} established the local and global (with small initial data) existence of strong solutions in some Besov spaces, respectively. In  the presence of  vacuum, under  the following  compatibility conditions,
\be\ba\label{tan}
\div u_0=\div H_0=0,\,\,\,-\Delta u_0+\na P_0-(H_0\cdot \na) H_0=\n_0^{1/2}g,~~~\mbox{in}~~\Omega,
\ea\ee
where $(P_0,~g)\in H^1\times L^2$ and $\Omega=\mathbb{R}^3$,
Chen-Tan-Wang \cite{tan} obtained the local existence of strong solutions to the 3D Cauchy problem, and proved the local solution is  global  provided the initial data satisfy some smallness conditions. When $\Omega \subset \mathbb{R}^2$ is a bounded domain,
 Huang-Wang \cite{hwjde1} investigated  the global existence
 of strong solution with general large data when the initial density contains vacuum states and the initial data satisfy the compatibility conditions \eqref{tan}.

Recently, Li-Liang \cite{LL2014} established the local existence of strong solutions to the 2D Cauchy problem of the compressible Navier-Stokes equations on the whole space $\rr$ with vacuum as far field density. Later,  L\"u-Huang \cite{lvhuang} obtained the local strong solutions to 2D Cauchy problem of the compressible MHD equations, which generalized the results of \cite{LL2014} to the MHD system. Motivated by \cite{LL2014}, Liang \cite{L2015} proved the local  existence of  strong solutions to 2D Cauchy problem of the incompressible Navier-Stokes equations, that is \eqref{mhd}-\eqref{n4} with $H=0$.
However, for the Cauchy problem  \eqref{mhd}-\eqref{n4} with $\Omega=\rr$, it is still open even for the local existence of strong  solutions when the far field   density is vacuum, in particular,  the initial density may have compact support. In fact, this is the main aim in this paper.

Now, we wish to define precisely what we mean by strong solutions.
\begin{definition}
If all derivatives involved in \eqref{mhd} for $(\rho,u,P,H)$ are regular distributions, and equations \eqref{mhd} hold almost everywhere in $\rr\times (0,T)$, then $(\n,u,P,H)$ is called a strong solution to \eqref{mhd}.
\end{definition}

In this section, for $1\le r\le \infty $ and $k\ge 1$, we denote the standard Lebesgue and Sobolev spaces as follows:
$$ L^r=L^r(\rr), \quad W^{k,r}= W^{k,r}(\rr), \quad H^k= W^{k,2}.  $$

\begin{theorem}\la{t1}
Let $\eta_0 $ be a positive constant and
\be\la{2.07}
\bar x\triangleq(e+|x|^2)^{1/2}\log^{1+\eta_0} (e+|x|^2).
\ee
For constants $q>2$ and $a>1$, assume that the initial data $(\n_0 ,u_0,H_0)$ satisfy
\be\ba\la{1.9}
\begin{cases}
\rho_0\geq0,\  \rho_0\bar x^a\in L^1 \cap H^1\cap W^{1,q},\,\,H_0\bar x^{a/2}\in L^2, \\
\sqrt{\rho_0}u_0\in L^2,\,\,\na u_0\in  L^2 , \,\,\na H_0\in  L^2 , \,\,\div u_0=\div H_0=0.
\end{cases}
\ea\ee
Then there exists a positive time $T_0>0$ such that the problem  \eqref{mhd}-\eqref{n4} has a unique strong solution $(\n,u,P,H)$ on $\rr\times (0,T_0]$ satisfying
\be\la{1.10}
\begin{cases}
0\le \rho\in C([0,T_0];L^1 \cap H^1\cap W^{1,q} ),\\
\rho\bar x^a\in L^\infty( 0,T_0 ;L^1\cap H^1\cap W^{1,q} ),\\
\sqrt{\n } u,\,\na u,\, \bar x^{-1}u,\,    \sqrt{t} \sqrt{\n}  u_t,\,    \sqrt{t} \na P,\,    \sqrt{t} \na^2  u \in L^\infty(0,T_0; L^2 ), \\
H,  H \bar{x}^{a/2}, \na H,  \sqrt{t}H_t,\,    \sqrt{t} \na^2 H \in L^\infty( 0,T_0 ;L^2), \\
\na u\in  L^2(0,T_0;H^1)\cap  L^{(q+1)/q}(0,T_0; W^{1,q}), \\
\na P\in  L^2(0,T_0;L^2)\cap  L^{(q+1)/q}(0,T_0;L^q), \\
\na H\in L^2(0, T_0; H^1),\,\,H_t,~\na H\bar{x}^{a/2}\in L^2(0, T_0; L^2),\\
\sqrt{t}\na u\in L^2(0,T_0; W^{1,q} ),  \\
\sqrt{\n} u_t, \,\sqrt{t}\na H\bar{x}^{a/2}, \,  \sqrt{t}\na u_t ,\, \sqrt{t}\na H_t ,\,  \sqrt{t} \bar x^{-1}u_t\in L^2(\rr\times(0,T_0)),\\
\end{cases}
\ee
and
\be \la{l1.2}
\inf\limits_{0\le t\le T_0}\int_{B_{N}}\n(x,t)dx\ge \frac14\int_{\rr} \n_0(x)dx,
\ee
for some constant $N >0$ and $B_{N }\triangleq\left.\left\{x\in\rr\right|
\,|x|<N \right\}$.
\end{theorem}

\begin{remark}\la{re4}
Compared  with \cite{tan} and \cite{CK2003}, there is no need to impose the additional compatibility conditions of the initial data for the local existence of strong solutions.
\end{remark}

If $H\equiv H_0\equiv0$, Theorem 1.1 directly yields the following local existence theorem for the density-dependent  Navier-Stokes equations.
\begin{theorem}\label{co1}
Let $\eta_0 $ and $\bar x$ be as in  \eqref{2.07}.
For constants $q>2$ and $a>1,$ assume that the initial data $(\n_0 ,u_0)$ satisfy
\bnn
\n_0\ge 0, \,  \rho_0\bar x^a\in   L^1 \cap H^1\cap W^{1,q},\,\,
\sqrt{\rho_0}u_0\in L^2,\,\,\na u_0\in  L^2 , \,\,\div u_0=0.
\enn
Then there exists a positive time $T_0>0$ such that  the 2D Cauchy  problem of the density-dependent  Navier-Stokes equations, thai is  \eqref{mhd}-\eqref{n4} with $H=0$,  has a unique strong solution $(\n,u,P)$ on $\rr\times (0,T_0]$ satisfying \eqref{1.10} where $H=0$, and \eqref{l1.2}.
\end{theorem}

\begin{remark} \la{re1}
Our Theorem \ref{co1} holds for arbitrary $a>1$ which is in sharp contrast to Liang \cite{L2015} where $a\in(1,2)$ is required.
\end{remark}

We now make some comments on the key ingredients of the analysis in this paper. It should be pointed out that, for the whole two-dimensional space, it seems difficult to bound the $L^p(\rr)$-norm of $u$ just in terms of
 $\|\n^{1/2} u\|_{L^2(\rr)} $ and   $\|\na u\|_{L^2(\rr)}$.
  Furthermore, as mentioned in many papers (see  \cite{hwjde1,lvhuang,lvshixu} for example),
 the strong coupling between the velocity field and the magnetic field, such as $\||u| |H|\|$ and  $\||u| |\na H|\|$,  will bring out some new difficulties. In order to overcome these difficulties stated above, we will use some key ideas due to \cite{LL2014,lvhuang} where they deal with the 2D compressible Navier-Stokes  and MHD equations, respectively.     On the one hand, motivated by \cite{LL2014}, it is enough to bound the  $L^p(\rr)$-norm of the momentum $\n u$ instead of  just the velocity $u$. More precisely, using a Hardy-type inequality (see \eqref{3.v2}) which is originally due to Lions  \cite{L1996}, together with some careful analysis on
 the spatial weighted estimate of the density (see \eqref{igj1-2}), we can obtain the desired  estimates on the $L^p(\rr)$-norm of $\n u$  (see \eqref{3.a2}). On the other hand, inspired by \cite{lvhuang}, we deduce  some spatial weighted estimates on both $H$ and $\na H$ (i.e., $\bar{x}^{a/2}H$ and $\bar{x}^{a/2}\na H$, see \eqref{lbqnew-gj10} and \eqref{igj10'}) which are crucial to control  the coupled  terms, such as $\||u| |H|\|$ and  $\||u| |\na H|\|$. Next, we then construct  approximate solutions to  \eqref{mhd}, that is, for density strictly away from vacuum initially,   consider   a initial boundary value problem of \eqref{mhd} in any bounded ball $B_R$  with radius $R>0.$  Finally, combining  all key points mentioned above with the similar arguments as in \cite{cho1,LL2014,lvhuang}, we derive some desired bounds on  the gradients of both the solutions and the spatial weighted  density,  which are independent of both the radius of the balls $B_R$ and the lower bound  of the initial density.

The rest of the paper is organized as follows: In Section 2, we collect some
elementary facts and inequalities which will be needed in later analysis. Sections 3 is devoted to the a priori estimates which are needed to obtain  the local existence and uniqueness of strong solutions. The main result
Theorem \ref{t1} is proved in Section 4.

\section{Preliminaries}

In this section, we will recall some  known facts and elementary
inequalities which will be used frequently later.
First of all, if the initial density is strictly away from vacuum, the following local existence theorem on bounded balls can be shown by similar arguments as in \cite{CK2003,tan,hwjde1}.
\begin{lemma}\la{th0}
For $R>0$ and $B_R=\{x\in\rr ||x|<R\}$, assume that $(\n_0,u_0,H_0 )$ satisfies
\be\ba\la{2.1}
&(\n_0,u_0,H_0) \in H^2(\O), \quad \inf\limits_{x\in\O}\n_0(x) >0 ,\quad\div u_0=\div H_0=0.
\ea\ee
Then there exist a small time $T_R>0$ such that the equations \eqref{mhd} with the following initial-boundary-value conditions
\be\la{ib2}
\begin{cases}
(\rho , u, H )(x,t=0)=(\n_0,u_0,H_0 ), \quad   &x\in \O,\\
u(x,t)=0,\, \, H(x,t)=0,\quad   &x\in\pa\O, \, t>0,
\end{cases}
\ee
has a unique classical solution $(\rho , u, P,H )$ on
$\O\times(0,T_R]$ satisfying
\be\la{lo1}
\begin{cases}
\n\in C\left([0,T_{R}];H^{2}\right),\\
(u, H) \in C\left([0,T_{R}]; H^{2}\right)\cap L^{2}\left(0,T_{R};H^{3}\right),\\
 P\in C\left([0,T_{R}]; H^{1}\right)\cap L^{2}\left(0,T_{R};H^{2}\right),
\end{cases}
\ee
where we denote $H^k=H^k(\O)$ for positive integer $k$.
\end{lemma}

Next, for $\Omega\subset\rr$, the following weighted $L^m$-bounds for elements of the Hilbert space $\ti D^{1,2}(\Omega)\triangleq\{v\in H^1_{\rm loc}(\Omega)|\na v\in L^2(\Omega)\}$ can be found in \cite[Theorem B.1]{L1996}.
\begin{lemma} \la{1leo}
For $m\in [2,\infty)$ and $\theta\in (1+m/2,\infty),$ there exists a positive constant $C$ such that for either $\Omega=\rr$ or $\Omega=B_R$ with $R\ge 1$ and for any $v\in \ti D^{1,2}(\Omega)$,
\be\la{3h}
\left(\int_{\Omega} \frac{|v|^m}{e+|x|^2}(\log (e+|x|^2))^{-\theta}dx  \right)^{1/m}\le C\|v\|_{L^2(B_1)}+C\|\na v\|_{L^2 (\Omega)}.
\ee
\end{lemma}

A useful consequence of Lemma \ref{1leo} is the following crucial weighted  bounds for elements of $\ti D^{1,2}(\Omega)$, which have been proved in  \cite[Lemma 2.4]{LL2014}.

\begin{lemma}\la{lemma2.6}
Let $\bar x$ and $\eta_0$ be as in \eqref{2.07} and $\Omega$ be as in Lemma \ref{1leo}. Assume that $\n \in L^1(\Omega)\cap L^\infty(\Omega)$ is a non-negative function such that
\be\la{2.i2}
\int_{B_{N_1} }\n dx\ge M_1,  \quad \|\n\|_{L^1(\Omega)\cap L^\infty(\Omega)}\le M_2,
\ee
for positive constants $M_1, M_2$, and $N_1\ge 1$ with $B_{N_1}\subset\Omega.$ Then  for $\ve> 0$ and $\eta>0,$ there is a positive constant $C$ depending only on $\ve,\eta, M_1,M_2,
N_1, $ and $\eta_0$ such that every $v\in \ti D^{1,2}(\Omega)$ satisfies
\be\ba\la{3.i2}
\|v\bar x^{-\eta}\|_{L^{(2+\ve)/\ti\eta}(\Omega)} &\le C \|\n^{1/2}v\|_{L^2(\Omega)}+C \|\na v\|_{L^2(\Omega)}
\ea\ee
with $\ti\eta=\min\{1,\eta\}$.
\end{lemma}

Finally, the following $L^p$-bound for elliptic systems, whose proof is similar to that of \cite[Lemma 12]{cho1}, is a direct result of the combination of the well-known elliptic theory \cite{adn} and a standard scaling procedure.
\begin{lemma} \la{lemma2.4}For $p>1$ and $k\ge 0,$ there exists a positive constant $C$ depending only on $p$ and $k$ such that
\be \la{lp} \|\na^{k+2}v\|_{L^p(\O)}\le C\|\Delta v\|_{W^{k,p}(\O)},\ee
for every $v\in W^{k+2,p}(\O)$ satisfying
\bnn v =0\ \ \mbox{on}\ \ \O.\enn
\end{lemma}

\section{A priori estimates}
In this section, for $r\in [1,\infty]$ and $k\ge0$, we denote
\bnn \int \cdot dx=\int_{\O}\cdot dx, \quad L^r=L^r(\O),\quad W^{k,r}=W^{k,r}(\O),\quad H^k=W^{k,2}.\enn
Moreover, for $R>4N_0\ge 4,$ assume that $(\n_0,u_0,H_0)$ satisfies, in addition to \eqref{2.1}, that
\be\la{w1}
1/2\le \int_{B_{N_0}}\n_0(x)dx \le \int_{B_R }\n_0(x)dx \le 3/2.
\ee
 Lemma \ref{th0} thus yields that  there exists some $T_R>0$ such that the  initial-boundary-value problem \eqref{mhd} and \eqref{ib2} has a unique classical solution $(\n,u,P,H)$ on $B_R\times[0,T_R]$ satisfying \eqref{lo1}.

Let $\bar{x}, \eta_0, a$, and $q$ be as in Theorem \ref{t1}, the main aim of this section is to derive the following key a priori estimate on $\psi $ defined by
\be\ba\la{3.2}
\psi(t)\triangleq & 1+\|\n^{1/2}u\|_{L^2}+\|\na u\|_{L^2}+\|\na H\|_{L^2} +\|\bar{x}^{a/2}H\|_{L^{2}}+\|\bar{x}^{a}\n\|_{L^{1}\cap H^{1}\cap W^{1,q}}.
\ea\ee
\begin{proposition} \la{pro}
Assume that $(\n_0,u_0,H_0)$ satisfies \eqref{2.1} and \eqref{w1}. Let $(\n,u,P,H)$ be the solution to the initial-boundary-value problem  \eqref{mhd} and \eqref{ib2}  on $\O\times (0,T_R]$ obtained by Lemma \ref{th0}. Then there exist positive constants $T_0$ and $M$ both depending only on $\mu,\nu, q$, $a$, $\eta_0$, $N_0,$ and $E_0$ such that
\be\la{o1}\ba
&\sup\limits_{0\le t\le T_0}\left(\psi(t)+\sqrt{t}\|\sqrt{\n}u_t\|_{L^2}+\sqrt{t}\|H_t\|_{L^2}+\sqrt{t}\|\na^2u\|_{L^2}+\sqrt{t}\|\na P\|_{L^2}+\sqrt{t}\|\na^2H\|_{L^2}\right)\\
&\quad+  \int_0^{T_0} \left( \|\sqrt{\n}u_t\|_{L^2}^2+\|\na^2u\|_{L^2}^{2}+ \|\na^2H\|_{L^2}^{2}+\|H_t\|_{L^2}^2+\|\na H\bar x^{a/2}\|_{L^2}^2\right) dt\\
&\quad+  \int_0^{T_0} \left(\norm[L^{q}]{\nabla^2u}^{(q+1)/q}+\norm[L^{q}]{\nabla P}^{(q+1)/q}+t\|\na^2u\|_{L^q}^{2} +t\|\na P\|_{L^q}^{2}\right)\\
&\quad+  \int_0^{T_0} \left( t\|\na u_t\|_{L^2}^{2} +t\|\na H_t\|_{L^2}^{2} \right) dt\le M,
\ea\ee
where
\bnn\ba E_0\triangleq &\|\n_0^{1/2}u_0\|_{L^2}+\|\na u_0\|_{L^2}+\|\na H_0\|_{L^2}+\|\bar{x}^{a}\n_0\|_{L^{1}\cap H^{1}\cap W^{1,q}}+ \|\bar{x}^{a/2}H_0\|_{L^2}.
\ea \enn
\end{proposition}

To show Proposition \ref{pro}, whose proof will be postponed to the end of this section, we begin with the following standard energy estimate for $(\n,u,P,H)$ and the estimate on the $L^p$-norm of the density.
\begin{lemma} \la{l3.00}
Under the conditions of Proposition \ref{pro}, let  $(\n,u,P,H)$ be a
smooth solution to  the  initial-boundary-value problem  \eqref{mhd} and \eqref{ib2}.  Then  for any $t>0$,
\be\ba\la{gj1}
\sup_{0\le s\le t} \left(\|\n\|_{L^1\cap L^\infty}+\|\rho^{1/2}u\|^2_{L^{2}}+\|H\|_{L^2}^2\right)
+\int_{0}^{t}\left(\|\na u\|_{L^2}^2+\|\na H\|_{L^2}^2 \right)ds \le C,
\ea\ee
where (and in what follows) $C$ denotes a generic positive constant depending only on $\mu,\nu, q, a$, $\eta_0$, $N_0,$ and $E_0$.
\end{lemma}

\pf First, applying standard energy estimate to \eqref{mhd} gives
\be\ba\la{gj1-1}
\sup_{0\le s\le t} \left(\|\rho^{1/2}u\|^2_{L^{2}}+\|H\|_{L^2}^2\right)
+\int_{0}^{t}\left(\|\na u\|_{L^2}^2+\|\na H\|_{L^2}^2 \right)ds  \le C.
\ea\ee
Next, since $\div u=0$, it is easy to deduce from  \eqref{mhd}$_1$  that (see \cite{L1996}),
\begin{equation}\label{3.2}
\sup_{0\le s\le t}\|\rho\|_{L^{1}\cap L^\infty}\leq C.
\end{equation}
This together with \eqref{gj1-1} yields \eqref{gj1} and completes the proof of Lemma \ref{l3.00}.   \hfill $\Box$

Next, we will give some spatial weighted estimates on the density and the magnetic.
\begin{lemma} \la{l3.01}
Under the conditions of Proposition \ref{pro}, let $(\n,u,P,H)$ be a
smooth solution to the initial-boundary-value problem \eqref{mhd} and \eqref{ib2}.
Then there exists a $T_1=T_1(N_0,E_0)>0$ such that for all $t\in (0, T_1],$
\be\ba\la{igj1}
\sup_{0\le s\le t} \left(\|\n\bar x^a\|_{L^1}+\|H\bar x^{a/2}\|_{L^2}^2\right)+\int_{0}^{t}\|\na H \bar x^{a/2}\|_{L^2}^2 ds\le C.
\ea\ee
\end{lemma}

\pf First, for $N>1$, let $\vp_N\in C^\infty_0(B_N)$ satisfy
\be\la{vp1}
0\le \vp_N \le 1, \quad \vp_N(x)=1, \mbox{ if } |x|\le N/2, \quad |\na \vp_N|\le C N^{-1}.
\ee
It follows from \eqref{mhd}$_1$ and  \eqref{gj1} that
\be\ba\la{oo0}
\frac{d}{dt}\int \n \vp_{2N_0} dx &=\int \n u \cdot\na \vp_{2N_0} dx\\
&\ge - C N_0^{-1}\left(\int\n dx\right)^{1/2}\left(\int\n |u|^2dx\right)^{1/2}\ge - \ti C (E_0).
\ea\ee
  Integrating
\eqref{oo0} and using \eqref{w1} give
\be\ba\la{p1}
\inf\limits_{0\le t\le T_1}\int_{B_{2N_0}}  \n dx&\ge \inf\limits_{0\le t\le T_1}\int \n \vp_{2N_0} dx\ge \int \n_0 \vp_{2N_0} dx-\ti C T_1  \ge 1/4.
\ea\ee
Here, $T_1\triangleq\min\{1, (4\ti C)^{-1}\} .$ From now on, we will always assume that $t\le T_1.$
The combination of \eqref{p1}, \eqref{gj1}, and \eqref{3.i2} implies that for $\ve> 0$ and $\eta>0$,  every $v\in \ti D^{1,2}(\O) $  satisfies\be \la{3.v2}\ba  \|v\bar x^{-\eta}\|_{L^{(2+\ve)/\ti\eta} }^2 &\le C(\ve,\eta) \|\n^{1/2}v\|_{L^2} +C (\ve,\eta)  \|\na v\|_{L^2 }^2, \ea\ee with $\ti\eta=\min\{1,\eta\}.$

 Next, multiplying $\eqref{mhd}_{1}$ by $\bar x^a$ and integrating by parts imply that
\be\ba\la{igj1-1}
\frac{d}{dt}\|\n  \bar x^a\|_{L^1}&\le C \int \n |u| \bar x^{a-1}\log^{1+\eta_0}(e+|x|^2)dx\\
&\le C\|\n \bar x^{a-1+\frac{8}{8+a}} \|_{L^\frac{8+a}{7+a}}\|u \bar x^{-\frac{4}{8+a}} \|_{L^{8+a}}\\
&\le C\|\n\|_{L^\infty}^{\frac{1}{8+a}}\|\n  \bar x^a\|_{L^1}^{\frac{7+a}{8+a}}\xl(\|\n^{1/2}u\|_{L^2} +\|\na u\|_{L^2 }\xr)\\
&\le C\xl(1+\|\n  \bar x^a\|_{L^1}\xr)\xl(1+\|\na u\|_{L^2 }^2\xr)
\ea\ee
due to \eqref{gj1} and \eqref{3.v2}. This combined   with Gronwall's inequality and \eqref{gj1} leads to \be\ba\la{igj1-2}
\sup_{0\le s\le t}\|\n  \bar x^a\|_{L^1}  \le  C\exp\left\{C\int_{0}^{t}\left(1+ \|\na u\|_{L^2}^2\right)ds\right\}\le C.
\ea\ee

Now, multiplying \eqref{mhd}$_3$ by $H\bar{x}^a$ and integrating by parts yield
\be\ba\la{lv4.1}
&\frac{1}{2}\frac{d}{dt}\|H\bar{x}^{a/2}\|_{L^2}^2+\nu \| \na H \bar{x}^{a/2}\|_{L^2}^2\\
&=\frac{\nu}{2}\int |H|^2\Delta\bar{x}^adx+\int H\cdot\na u\cdot H\bar{x}^adx +\frac{1}{2}\int |H|^2u\cdot\na\bar{x}^adx\\
&\triangleq \bar{I}_1+\bar{I}_2+\bar{I}_3,
\ea\ee
where
\be\ba\la{lv4.4}
|\bar{I}_1| &\leq C\int |H|^2 \bar{x}^a \bar{x}^{-2}\log^{2(1-\eta_0)}(e+|x|^2) dx\leq C\int |H|^2 \bar{x}^a dx,\\
|\bar{I}_2|&
\leq C \|\na u\|_{L^2}\| H   \bar{x}^{a/2} \|_{L^4}^2\\
&\leq   C \|\na u\|_{L^2}\| H   \bar{x}^{a/2}\|_{L^2} (\| \na H \bar{x}^{a/2} \|_{L^2}+\| H \na\bar{x}^{a/2}\|_{L^2})\\
&\leq C(\|\na u\|_{L^2}^2+1)\| H \bar{x}^{a/2}\|_{L^2}^2 +\frac\nu4\| \na H \bar{x}^{a/2}\|_{L^2}^2,\\
|\bar{I}_3|
&\leq    C\| H \bar{x}^{a/2}\|_{L^4}\| H \bar{x}^{a/2}\|_{L^2} \|u\bar{x}^{-3/4}\|_{L^{4}}\\
&\leq   C\| H \bar{x}^{a/2}\|_{L^4}^2+C\|H \bar{x}^{a/2}\|_{L^2}^2 \left( \|\n^{1/2}u\|_{L^2}^2+\|\na u\|_{L^2}^2\right) \\
&\leq C\left(1+ \|\na u\|_{L^2}^2\right)\|H\bar{x}^{a/2}\|_{L^2}^2+\frac\nu4\| \na H \bar{x}^{a/2}\|_{L^2}^2,
\ea\ee
due to Gagliardo-Nirenberg inequality (see \cite{nir}), \eqref{gj1}, and \eqref{3.v2}. Putting \eqref{lv4.4} into \eqref{lv4.1}, we get after  using Gronwall's inequality and \eqref{gj1} that
\be\ba\la{lbqnew-gj10}
\sup_{0\le s\le t} \|H\bar{x}^{a/2}\|_{L^2}^2 +\int_{0}^{t}\| \na H \bar{x}^{a/2}\|_{L^2}^2ds\le C\exp\left\{C\int_{0}^{t}\left(1+ \|\na u\|_{L^2}^2\right)ds\right\}\le C,
\ea\ee
which together with \eqref{igj1-2} gives  \eqref{igj1} and finishes the proof of Lemma \ref{l3.01}.   \hfill $\Box$

\begin{lemma} \la{l3.0}  Let $(\n,u,P,H)$ and $T_1$ be as in Lemma \ref{l3.01}. Then there exists a positive constant $\alpha>1$ such that for all $t\in (0, T_1]$,
\be\ba\la{gj3}
&\sup_{0\le s\le t}\left(\|\na u\|_{L^2}^2+\|\na H\|_{L^2}^2\right)\\
&\quad+\int_{0}^{t}\left(\|\n^{1/2}u_s\|_{L^2}^2+\|\na^2 u\|_{L^2}^2+\| H_s\|_{L^2}^2+\|\na^2 H\|_{L^2}^2\right)ds\le C+C\int_{0}^{t}\psi^\alpha(s)ds.
\ea\ee
\end{lemma}

\pf Multiplying $\eqref{mhd}_{2}$ by $u_t$ and integrating by parts, one has
\be\ba\la{3r1}
& \mu\frac{d}{dt}\int|\na u|^2dx+\int \n|  u_t|^2dx \le C\int\n|u|^2|\na u|^2dx+\int H\cdot \na H \cdot u_t dx.
\ea\ee

First,  it follows from \eqref{gj1}, \eqref{igj1}, and \eqref{3.v2} that for any  $\ve> 0$ and any $\eta>0 ,$
\be \la{local1}\ba \|\n^\eta v\|_{L^{(2+\ve)/\ti\eta}}
&\le C\|\n^\eta \bar x^{\frac{3\ti\eta a}{4(2+\ve)}} \|_{L^{ \frac{4(2+\ve)}{3\ti\eta}}} \|v\bar x^{-\frac{3\ti\eta a}{4(2+\ve)}} \|_{L^{ \frac{4(2+\ve)}{\ti\eta}}} \\
&\le C\left(\int \n^{\frac{4(2+\ve)\eta}{3\ti\eta}-1}\n \bar x^a dx\right)^{ \frac{3\ti\eta}{4(2+\ve)}} \|v\bar x^{-\frac{3\ti\eta a}{4(2+\ve)}} \|_{L^{ \frac{4(2+\ve)}{\ti\eta}}} \\
&\le C\|\n\|_{L^\infty}^{\frac{4(2+\ve)\eta-3\ti\eta}{4(2+\ve)}}
\|\n\bar x^a\|_{L^1}^{\frac{3\ti\eta}{4(2+\ve)}}
\left( \|\n^{1/2}v\|_{L^2 } +\|\na v\|_{L^2 }\right) \\
&\le C\|\n^{1/2}v\|_{L^2 } +C\|\na v\|_{L^2 },
\ea\ee
where $\ti\eta=\min\{1,\eta\}$ and $v\in \ti D^{1,2}(\O) $. In particular, this together with \eqref{gj1} and \eqref{3.v2} derives
\be \la{3.a2}\ba&\|\n^\eta u\|_{L^{(2+\ve)/\ti\eta}}+ \|u\bar x^{-\eta}\|_{L^{(2+\ve)/\ti\eta}} \le C(1+\|\na u\|_{L^2 } ),
\ea\ee
which combined with H\"older's and Gagliardo-Nirenberg inequalities yields
\be\ba\la{cc5}
\int\n |u|^2|\na u|^2dx &\le C\|\n^{1/2}u\|_{L^8}^{2}\|\na u\|_{L^{8/3}}^{2} \\ &\le C\|\n^{1/2}u\|_{L^8}^{2}\|\na u\|_{L^2}^{3/2} \|\na  u\|_{H^1}^{1/2}\\
&\le C \psi^{\alpha }+\ve \|\na^2 u\|_{L^2}^2,
\ea\ee
where (and in what follows) we use $\alpha>1$ to denote a genetic constant, which may be different from line to line.

For the second term on the right-hand side of \eqref{3r1}, integration by parts together with \eqref{mhd}$_4$ and Gagliardo-Nirenberg inequality deduces   that
\be\ba\la{bv11}
\int H\cdot \na H \cdot u_tdx 
&=-\frac{d}{dt} \int H\cdot \na u \cdot Hdx +\int H_t\cdot \na u\cdot H dx+\int H\cdot\na u\cdot H_tdx\\
&\le -\frac{d}{dt} \int H\cdot \na u \cdot Hdx+ \frac{\nu^{-1}}{2} \|H_t\|_{L^2}^2+C \|H\|_{L^4}^2\|\na u\|_{L^4}^2\\
&\le -\frac{d}{dt} \int H\cdot \na u \cdot Hdx+\frac{\nu^{-1}}{2}\|H_t\|_{L^2}^2+C \|H\|_{L^2} \|\na H\|_{L^2}\|\na u\|_{L^2}\|\na   u\|_{H^1}\\
&\le -\frac{d}{dt} \int H\cdot \na u \cdot Hdx+\frac{\nu^{-1}}{2} \|H_t\|_{L^2}^2+\ve \|\na^2 u\|_{L^2}^2+C \psi^{\alpha}.
\ea\ee

Inserting \eqref{cc5}--\eqref{bv11} into \eqref{3r1} gives
\be\la{3.8-1}\ba & \frac{d}{dt}B(t)+\int \n|  u_t|^2dx \le \ve \|\na^2 u\|_{L^2}^2+\frac{\nu^{-1}}{2} \|H_t\|_{L^2}^2+C  \psi^\alpha,\ea\ee
where
$$B(t)\triangleq\mu\|\na u\|_{L^2}^2+\int H\cdot \na u \cdot Hdx$$
satisfies
\be\ba\la{lv4.8'}
&\frac{\mu}{2}\|\na u\|_{L^2}^2-C_1\|\na H\|_{L^2}^2 \le B(t) \leq  C\|\na u\|_{L^2}^2+ C\|\na H\|_{L^2}^2.
\ea\ee

Moreover, it follows from \eqref{mhd}$_3$  that
\be\la{lv4.9}\ba
& \nu\frac{d}{dt}\|\na H\|_{L^2}^2 +\|H_t\|_{L^2}^2+ \nu^{2}\|\Delta H\|_{L^2}^2 \\
& \le C\||H| |\na u|\|_{L^2}^2+C\||u| |\na H|\|_{L^2}^2 \\
&\le C\|H\|_{L^2} \|\na^2 H\|_{L^2}\|\na u\|_{L^2}^2+C\|\bar{x}^{-a/4}u\|_{L^8}^2\|\bar{x}^{a/2}\na H\|_{L^2} \|\na H\|_{L^4}\\
&\le  \frac{\nu^2}{2}\|\Delta H\|_{L^2}^2+C\psi^\alpha+C\|\bar{x}^{a/2}\na H\|_{L^2}^2
\ea\ee
due to \eqref{lp}, \eqref{3.a2}, and Gagliardo-Nirenberg inequality. Multiplying \eqref{lv4.9} by $\nu^{-1}(C_1+1)$  and adding the resulting inequality to \eqref{3.8-1} imply
\be\ba\la{lv4.10}
& \frac{d}{dt} \left(B(t)+(C_1+1)\|\na H\|_{L^2}^2\right) +\|\n^{1/2}u_t\|_{L^2}^2+\frac{\nu^{-1}}{2}\|H_t\|_{L^2}^2+ \frac{\nu }{2} \|\Delta H\|_{L^2}^2\\
&\le C\psi^\alpha+C\|\bar{x}^{a/2}\na H\|_{L^2}^2+\ve \|\na^2 u\|_{L^2}^2.
\ea\ee

On the other hand, since $(\n,u,P, H)$ satisfies the following Stokes system
\be\la{stokes1}
\begin{cases}
-\mu\Delta u + \nabla P = -\n u_t-\n u\cdot\na u+ H\cdot\na H-\frac{1}{2}\na |H|^2,\,\,\,\,&x\in B_R,\\
\div u=0,   \,\,\,&x\in  B_R,\\
u(x)=0,\,\,\,\,&x\in \p B_R,
\end{cases}
\ee
applying the standard $L^p$-estimate to \eqref{stokes1} (see \cite{T2001}) yields that for any $p\in (1,\infty),$
\be\ba\label{stokes2}
\|\na^2 u \|_{L^p}+\|\nabla P\|_{L^p}\le C\|\rho u_t\|_{L^p}+C\|\rho u\cdot\na u\|_{L^p}+ C\||H||\na H|\|_{L^p}.
\ea\ee
Then, it follows from \eqref{stokes2}, \eqref{gj1}, \eqref{3.a2}, and Gagliardo-Nirenberg inequality that
\be\ba\label{001}
&\|\na^2 u \|_{L^2}^2+\|\nabla P\|_{L^2}^2\\
&\leq C\|\rho u_t\|_{L^2}^2+C\|\rho u\cdot\na u\|_{L^2}^2+ C\||H||\na H|\|_{L^2}^2\\
&\leq  C\|\rho\|_{L^\infty} \|\sqrt{\rho} u_t\|_{L^2}^2+C\|\rho u\|_{L^4}^2\|\na u\|_{L^4}^2+ C\|H\|_{L^4}^2\|\na H\|_{L^4}^2\\
&\leq  C \|\sqrt{\rho} u_t\|_{L^2}^2+C\|\rho u\|_{L^4}^2\|\na u\|_{L^2} \|\na u\|_{H^1} + C\|H\|_{L^2} \|\na H\|_{L^2}^2\|\na H\|_{H^1} \\
&\leq C \|\sqrt{\rho} u_t\|_{L^2}^2+\frac{1}{4}\|\na^2 H\|_{L^2}^2+\frac{1}{2}\|\na^2 u\|_{L^2}^2+C\left(1+\|\na H\|_{L^2}^4+\|\na u\|_{L^2}^6\right)\\
&\leq C \|\sqrt{\rho} u_t\|_{L^2}^2+\frac{1}{4}\|\na^2 H\|_{L^2}^2+\frac{1}{2}\|\na^2 u\|_{L^2}^2+ C\psi^\alpha.
\ea\ee

Finally, substituting \eqref{001} into \eqref{lv4.10}  and choosing $\ve$ suitably small, one gets
\be\ba\la{lvzh4.10}
& \frac{d}{dt} \left(B(t)+(C_1+1)\|\na H\|_{L^2}^2\right) +\frac{1}{2}\|\n^{1/2}u_t\|_{L^2}^2+\frac{\nu^{-1}}{2}\|H_t\|_{L^2}^2+ \frac{\nu }{4} \|\Delta H\|_{L^2}^2\\
&\le C\psi^\alpha+C\|\bar{x}^{a/2}\na H\|_{L^2}^2.
\ea\ee
Integrating the above inequality over $(0,t)$, it follows from \eqref{lp}, \eqref{lv4.8'}, \eqref{igj1}, and \eqref{001} that \eqref{gj3} holds. The proof of Lemma \ref{l3.0} is finished.  \hfill $\Box$

\begin{lemma} \la{le-3}
Let $(\n,u,P,H)$ and $T_1$ be as in Lemma \ref{l3.01}. Then there exists a positive constant $\alpha>1$ such that for all $t\in (0, T_1],$
\be\ba\la{gj6}
&\sup_{0\le s\le t}  \left(s\| \n^{1/2}u_s \|_{L^2}^2+s\|H_s\|_{L^2}^2\right) +\int_{0}^{t} \left( s\| \na u_s\|_{L^2}^2+s\| \na H_s\|_{L^2}^2\right) ds\\
&\le C\exp\left\{C\int_{0}^{t} \psi^\alpha ds\right\}.
\ea\ee
\end{lemma}

{\it Proof.} Differentiating $\eqref{mhd}_2$ with respect to $t$ gives
\be\ba\la{zb1}
&\n u_{tt}+\n u\cdot \na u_t-\mu\Delta u_t  \\
&=-\n_t(u_t+u\cdot\na u)-\n u_t\cdot\na u -\na P_t+\left(H\cdot\na H-\frac{1}{2}\na|H|^2\right)_t.
\ea\ee
Multiplying \eqref{zb1} by $u_t$ and integrating the resulting equality by parts over $\O,$ we obtain after using $\eqref{mhd}_1$ and $\eqref{mhd}_4$ that
\be\ba \la{na8}
&\frac{1}{2}\frac{d}{dt} \int \n |u_t|^2dx+\mu\int  |\na u_t|^2 dx \\
&\le C\int  \n |u||u_{t}| \left(|\na  u_t|+|\na u|^{2}+|u||\na^{2}u|\right)dx +C\int \n |u|^{2}|\na u ||\na u_{t}|dx \\
&\quad+C\int \n |u_t|^{2}|\na u |dx +\int H_t\cdot \na H\cdot u_tdx+\int H\cdot \na H_t\cdot u_t dx \\
&\triangleq \sum_{i=1}^5 \hat{I}_i.
\ea\ee

We estimate each term  on the right-hand side of  \eqref{na8} as follows.

First, it follows from   \eqref {local1}, \eqref{3.a2}, and  Gagliardo-Nirenberg inequality  that
\be\ba\la{na2}
\hat{I}_1& \le C \|\n^{1/2} u\|_{L^{6}}\|\n^{1/2} u_{t}\|_{L^{2}}^{1/2} \|\n^{1/2} u_{t}\|_{L^{6}}^{1/2}\left(\| \na u_{t}\|_{L^{2}}+\| \na u\|_{L^{4}}^{2} \right) \\
&\quad +C\|\n^{1/4}  u \|_{L^{12}}^{2}\|\n^{1/2} u_{t}\|_{L^{2}}^{1/2} \|\n^{1/2} u_{t}\|_{L^{6}}^{1/2} \| \na^{2} u \|_{L^{2}}  \\
& \le C(1+\| \na u \|_{L^{2}}^2) \|\n^{1/2} u_{t}\|_{L^{2}}^{1/2}\left(\|\n^{1/2} u_{t}\|_{L^{2}} +\| \na u_{t}\|_{L^{2}}\right)^{1/2}\\
&\qquad\cdot\left(\| \na u_{t}\|_{L^{2}}+  \| \na u \|^{2}_{L^{2}}+\|\na u \|_{L^{2}}\| \na^{2} u \|_{L^{2}}+ \| \na^{2} u \|_{L^{2}}\right)\\
&\le  \frac{\mu}{6}\| \na u_{t}\|_{L^{2}}^{2}+C\psi^{\alpha}\|\n^{1/2} u_{t}\|_{L^{2}}^{2}+C\psi^{\alpha}+ C\left(1+\| \na u \|_{L^{2}}^2\right)\| \na^{2} u \|_{L^{2}}^2.
\ea\ee
Then, H\"older's inequality combined with \eqref{local1} and \eqref{3.a2} leads to
\be\ba\la{5.a3}
\hat{I}_2+\hat{I}_3&\le C \|\n^{1/2} u\|_{L^{8}}^{2}\|\na u\|_{L^{4}} \| \na u_{t}\|_{L^{2}}+C \| \na u\|_{L^{2}}
\|\n^{1/2} u_{t}\|_{L^{6}}^{3/2}\|\n^{1/2} u_{t}\|_{L^{2}}^{1/2} \\
&\le \frac{\mu}{6} \| \na u_{t}\|_{L^{2}}^{2} + C \psi^{\alpha} \|\n^{1/2} u_{t}\|_{L^{2}}^{2}+ C \left(\psi^{\alpha}+\|\na^{2} u \|_{L^{2}}^{2}\right).
\ea\ee
Next,  integration by parts together with \eqref{mhd}$_4$, H\"older's and  Gagliardo-Nirenberg inequalities  indicates that
\be\ba\la{na3}
\hat{I}_4+\hat{I}_5 & =-\int H_t\cdot \na u_t\cdot Hdx-\int H\cdot \na u_t\cdot H_tdx \\
&\le \frac{\mu}{6} \|\na u_{t}\|_{L^{2}}^{2}+ C \|H\|_{L^{4}}^{2}\|H_t\|_{L^{4}}^{2}\\
&\le \frac{\mu}{6}\|\na u_{t}\|_{L^{2}}^{2}+\frac{\mu\nu}{4(C_2+1)} \|\na H_{t}\|_{L^{2}}^{2}+ C \psi^\alpha\|H_t\|_{L^{2}}^2,
\ea\ee
where the constant $C_2$ is defined in the following \eqref{ilv4.14}.

Substituting \eqref{na2}--\eqref{na3} into \eqref{na8}, we obtain after using   \eqref{001} that
\be\ba\la{a4.6}
\frac{d}{dt} \| \n^{1/2}u_{t}\|_{L^{2}}^{2}+\mu\|\na u_{t}\|_{L^{2}}^{2} & \le   C\psi^{\alpha}\left(1+\| \n^{1/2}u_{t}\|_{L^{2}}^{2}+\| H_t\|_{L^{2}}^{2} \right)\\
&\quad+  \frac{\mu\nu}{2(C_2+1)}\|\na H_t\|_{L^{2}}^{2}+C\xl(1+\|\na u\|_{L^{2}}^2\xr)\|\na^2 H\|_{L^{2}}^2.
\ea\ee

Differentiating $\eqref{mhd}_3$ with respect to $t$ shows
\be\la{lv4.12}\ba
H_{tt}-H_t\cdot\na u-H\cdot\na u_t+u_t\cdot\na H+u\cdot\na H_t=\nu\Delta H_t.
\ea\ee
Multiplying \eqref{lv4.12} by $H_t$ and integrating the resulting equality  over $\O$ yield that
\be\la{lv4.13}\ba
&\frac{1}{2}\frac{d}{dt} \int |H_t|^2dx+\nu \int  |\na H_t|^2dx\\
&=\int H\cdot\na u_t\cdot H_tdx-\int u_t\cdot\na H\cdot H_tdx +\int  H_t\cdot\na u\cdot H_tdx-\int u\cdot\na H_t\cdot H_tdx \\
&\triangleq \sum_{i=1}^{4} S_i.
\ea\ee

On the one hand, we deduce from \eqref{3.v2} and \eqref{lbqnew-gj10} that
\be\ba\la{lvbo4.14}
\sum_{i=1}^{2}  S_i
&\le C\|\na u_t\|_{L^2} \|H_t\|_{L^4} \|H\|_{L^4}+ C\|\na H_t\|_{L^2} \||u_t||H|\|_{L^2} \\
&\le C\|H_t\|_{L^4}^2+C\|\na u_t\|_{L^2}^2+ \frac{\nu}{8} \|\na H_t\|_{L^2}^2+C \||u_t||H|\|_{L^2}^2 \\
&\le \frac{\nu}{4} \|\na H_t\|_{L^2}^2+C \|H_t\|_{L^2}^2+C\|\na u_t\|_{L^2}^2+ C \|u_t \bar{x}^{-a/4}\|_{L^8}^2 \|H \bar{x}^{a/2}\|_{L^2}\|H\|_{L^4}\\
&\le \frac{\nu}{4} \|\na H_t\|_{L^2}^2 +C \|H_t\|_{L^2}^2+C  \|\na u_t\|_{L^2}^2+ C  \|\n^{1/2}u_t\|_{L^2}^2,
\ea\ee
where one has used the following estimate
\be\ba\la{gj2}
\sup_{0\le s\le t} \||H|^2\|_{L^2}^2
+\int_{0}^{t}\||\na H||H|\|_{L^2}^2ds \le C.
\ea\ee
Indeed, multiplying $\eqref{mhd}_{3}$ by $H|H|^2$ and integrating by parts lead to
\be\ba\la{lv3.7}
&\frac{1}{4}\left(\||H|^2\|^2_{L^2}\right)_t +\nu\||\na H| |H|\|^2_{L^2}+\frac{\nu}{2}\|\na |H|^2 \|^2_{L^2}\\
&\leq C\|\na u\|_{L^2} \||H|^2\|_{L^4}^2 \leq C\|\na u\|_{L^2} \||H|^2\|_{L^2} \|\na |H|^2\|_{L^2}\\
&\leq \frac{\nu}{4} \|\na |H|^2\|_{L^2}^2+ C\|\na u\|_{L^2}^2\||H|^2\|_{L^2}^2,
\ea\ee
which together with Gronwall's inequality and $\eqref{gj1}$ gives \eqref{gj2}.

On the other hand, integration by parts combined with $\eqref{mhd}_{4}$ and  Gagliardo-Nirenberg inequality yields
\be\ba\la{lv4.14}
\sum_{i=3}^{4}  S_i &=\int  H_t\cdot\na u\cdot H_tdx
\leq  C \|H_t\|_{L^2} \|\na H_t\|_{L^2} \|\na u\|_{L^2}
\leq \frac{\nu}{4}\|\na H_t\|_{L^2}^2+C\psi^{\alpha}\|H_t\|_{L^2}^2.
\ea\ee
Inserting \eqref{lvbo4.14} and \eqref{lv4.14} into \eqref{lv4.13}, one has
\be\la{ilv4.14}\ba
&\frac{d}{dt} \|H_t\|_{L^2}^2+\nu\|\na H_t\|_{L^2}^2\le  C\psi^{\alpha}\xl(\|H_t\|_{L^2}^2+\|\n^{1/2}u_t\|_{L^2}^2\xr)+C_2  \|\na u_t\|_{L^2}^2.
\ea\ee

Finally, multiplying \eqref{a4.6} by $\mu^{-1}(C_2+1) $ and adding the resulting inequality to \eqref{ilv4.14}, we get
\be\ba\la{ilv4.13}
&\frac{d}{dt}  \left(\mu^{-1}(C_2+1)\|\n^{1/2}u_t\|_{L^2}^2+\|H_t\|_{L^2}^2\right)+\|\na u_t\|_{L^2}^2+\frac{\nu}{2}\|\na H_t\|_{L^2}^2\\
&\le  C\psi^{\alpha}\xl(1+\|H_t\|_{L^2}^2+  \|\n^{1/2}u_t\|_{L^2}^2\xr) +C\xl(1+\|\na u\|_{L^{2}}^2\xr)\|\na^2 H\|_{L^{2}}^2.
\ea\ee
Multiplying \eqref{ilv4.13} by $t$, we obtain \eqref{gj6} after using  Gronwall's inequality and \eqref{gj3}. The proof of Lemma \ref{le-3} is finished.  \hfill $\Box$


\begin{lemma} \la{le-3'}
Let $(\n,u,P,H)$ and $T_1$ be as in Lemma \ref{l3.01}. Then there exists a positive constant $\alpha>1$ such that for all $t\in(0, T_1],$
\be\ba\la{gj10'}
&\sup_{0\le s\le t} \left(s\|\na^2 u\|_{L^2}^2+s\|\na^2 H\|_{L^2}^2+s\|\na H\bar{x}^{a/2}\|_{L^2}^2\right)+\int_{0}^{t} s\|\Delta H\bar{x}^{a/2}\|_{L^2}^2ds\\
&\le C\exp\left\{C\exp\left\{C\int_{0}^{t} \psi^{\alpha} ds\right\}\right\}.
\ea\ee
\end{lemma}

\pf First, multiplying \eqref{mhd}$_3$ by $\Delta H\bar{x}^a$ and  integrating by parts lead to
\be\ba\la{AMSS5}
&\frac{1}{2}\frac{d}{dt}\int |\na H|^2\bar{x}^adx+\nu \int |\Delta H|^2\bar{x}^adx\\
\le& C\int|\na H| |H| |\na u| |\na\bar{x}^a|dx+C\int|\na H|^2|u| |\na\bar{x}^a|dx+C\int|\na H| |\Delta H| |\na \bar{x}^a|dx\\
&+C\int |H||\na u||\Delta H|\bar{x}^adx+C\int |\na u||\na H|^2  \bar{x}^adx
\triangleq \sum_{i=1}^5 J_i.
\ea\ee

Using \eqref{lbqnew-gj10}, \eqref{3.v2}, H\"older's and Gagliardo-Nirenberg inequalities, one gets by some direct calculations that
\be\ba
J_1\le & C\|H\bar{x}^{a/2}\|_{L^4}\|\na u\|_{L^4}\|\na H\bar{x}^{a/2}\|_{L^2}\\
\le & C\|H\bar{x}^{a/2}\|_{L^2}^{1/2}\left(\|\na H\bar{x}^{a/2}\|_{L^2}+\|H\bar{x}^{a/2}\|_{L^2}\right)^{1/2}\|\na u\|_{L^2}^{1/2}\|\na u\|_{H^1}^{1/2}\|\na H\bar{x}^{a/2}\|_{L^2}\\
\le &C\psi^{\alpha}+C\|\na^2 u\|_{L^2}^2+C\psi^{\alpha}\|\na H\bar{x}^{a/2}\|_{L^2}^2,\\
J_2\leq & C\||\na H|^{2-\frac{2}{3a}}\bar{x}^{a-\frac{1}{3}}\|_{L^{\frac{6a}{6a-2}}} \|u\bar{x}^{-\frac{1}{3}}\|_{L^{6a}}\||\na H|^{\frac{2}{3a} }\|_{L^{6a}}\\
\le &C\psi^{\alpha}\|\na H \bar{x}^{a/2} \|_{L^2}^\frac{6a-2}{3a}\|\na H \|_{L^4}^\frac{2}{3a}  \leq  C\psi^{\alpha}\|\na H \bar{x}^{a/2} \|_{L^2}^2+  C\|\na H \|_{L^4}^2\\
\leq  &C\psi^{\alpha}\|\na H \bar{x}^{a/2}\|_{L^2}^2+   \frac{\nu}{4}\|\Delta H \bar{x}^{a/2}\|_{L^2}^2,\\
J_3+J_4\le &\frac{\nu}{4}\|\Delta H\bar{x}^{a/2}\|_{L^2}^2+C\|\na H\bar{x}^{a/2}\|_{L^2}^2
+C\|H\bar{x}^{a/2}\|_{L^4}^2 \|\na u\|_{L^4}^2\\
\le &\frac{\nu}{4}\|\Delta H\bar{x}^{a/2}\|_{L^2}^2+C\|\na H\bar{x}^{a/2}\|_{L^2}^2\\
& +C\|H\bar{x}^{a/2}\|_{L^2}\left(\|\na H\bar{x}^{a/2}\|_{L^2}+\|H\bar{x}^{a/2}\|_{L^2}\right)\|\na u\|_{L^2} \|\na u\|_{H^1} \\
\le & \varepsilon\|\Delta H\bar{x}^{a/2}\|_{L^2}^2+C\psi^{\alpha}\|\na H\bar{x}^{a/2}\|_{L^2}^2+C\psi^{\alpha}+C\|\na^2 u\|_{L^2}^2,\\
J_5\le &C\|\na u\|_{L^\infty} \|\na H\bar{x}^{a/2}\|_{L^2}^2
\le  C\left(\psi^{\alpha}+\|\na^2 u\|_{L^q}^{(q+1)/q}\right)\|\na H\bar{x}^{a/2}\|_{L^2}^2.
\ea\ee
Substituting the above estimates   into \eqref{AMSS5} gives
\be\ba\la{AMSS10}
& \frac{d}{dt}\int |\na H|^2\bar{x}^adx+\nu \int |\Delta H|^2\bar{x}^adx\\
&\quad\le C\left(\psi^{\alpha}+\|\na^2 u\|_{L^q}^{(q+1)/q}\right)\|\na H\bar{x}^{a/2}\|_{L^2}^2+C\left(\|\na^2 u\|_{L^2}^2+\psi^{\alpha}\right).
\ea\ee

Now, we claim that
\be\ba\la{gj7}
&\int_{0}^{t} \left( \| \na^2 u\|_{ L^q}^{(q+1)/q} + \| \na P\|_{ L^q}^{(q+1)/q} + s\| \na^2 u\|_{ L^q}^2+ s\| \na P\|_{ L^q}^2\right)  ds\\
&\le C\exp\left\{C\int_{0}^{t} \psi^\alpha(s)ds\right\},
\ea\ee
whose proof will be given at the end of this proof. Thus, multiplying \eqref{AMSS10} by $t$,  we infer from \eqref{lbqnew-gj10}, \eqref{gj3}, \eqref{gj7}, and Gronwall's inequality  that
\be\ba\la{igj10'}
&\sup_{0\le s\le t} \left( s\|\na H\bar{x}^{a/2}\|_{L^2}^2\right)+\int_{0}^{t} s\|\Delta H\bar{x}^{a/2}\|_{L^2}^2ds \le C\exp\left\{C\exp\left\{C\int_{0}^{t} \psi^{\alpha} ds\right\}\right\}.
\ea\ee

Next,  it is easy to deduce from $\eqref{mhd}_3$, \eqref{lp}, \eqref{gj1-1}, \eqref{3.a2}, H\"older's and Gagliardo-Nirenberg inequalities that
\be\ba\la{AMSS11}\|\na^2 H\|^2_{L^2}&\leq C\|H_t\|^2_{L^2}+C\||u| |\na H|\|_{L^2}^2+C\||H| |\na u|\|^2_{L^2}\\
&\leq C\|H_t\|^2_{L^2}+C\|u \bar{x}^{-a/4}\|_{L^8}^2\|\na H \bar{x}^{a/2}\|_{L^2}\|\na H\|_{L^4}+C\|H\|_{L^2}\|\na^2 H\|_{L^2}\|\na u\|^2_{L^2}\\
&\leq C\|H_t\|^2_{L^2}+C\|\na H \bar{x}^{a/2}\|_{L^2}^2+C\|u \bar{x}^{-a/4}\|_{L^8}^4\|\na H\|_{L^4}^2+C\|\na^2 H\|_{L^2} \|\na u\|^2_{L^2}\\
&\leq C\|H_t\|^2_{L^2}+C\|\na H \bar{x}^{a/2}\|_{L^2}^2+\frac{1}{4}\|\na^2 H\|_{L^2}^2 +C\xl(1+\|\na u\|^8_{L^2}\xr)\xl(1+\|\na H\|^2_{L^2}\xr),
\ea\ee
which together with \eqref{001} gives that
\be\ba\la{iAMSS11}
&\|\na^2u \|_{L^2}^2 +\|\na P \|_{L^2}^2+\|\na^2 H\|^2_{L^2} \\
&\leq C\xl(\|\sqrt{\rho} u_t\|_{L^2}^2+\|H_t\|^2_{L^2}+ \|\na H \bar{x}^{a/2}\|_{L^2}^2\xr) +C\xl(1+\|\na u\|^8_{L^2}\xr)\xl(1+\|\na H\|^4_{L^2}\xr).
\ea\ee
Then, multiplying \eqref{iAMSS11} by $s$, one gets from  \eqref{gj3}, \eqref{gj6}, and \eqref{igj10'} that
\be\ba\la{iAMSS12}
&\sup_{0\le s\le t} \left( s\|\na^2u \|_{L^2}^2 + s\|\na P \|_{L^2}^2 +s\|\na^2 H\|^2_{L^2} \right)\\
&\leq C\exp\left\{C\exp\left\{C\int_{0}^{t} \psi^{\alpha} ds\right\}\right\}  +C\left(1+\int_{0}^{t} \psi^\alpha(s)ds\right)^{12}\\
&\leq C\exp\left\{C\exp\left\{C\int_{0}^{t} \psi^{\alpha} ds\right\}\right\},
\ea\ee
which combined with \eqref{igj10'} implies \eqref{gj10'}.

Finally, to finish the proof of Lemma \ref{le-3'}, it suffices to show \eqref{gj7}. Indeed, choosing $p=q$ in \eqref{stokes2}, we deduce from    \eqref{gj1}, \eqref{local1},  and Gagliardo-Nirenberg inequality that
\be\ba\la{lv3.60}
 \|\na^2u\|_{L^q}+\|\na P\|_{L^q}
& \le C \left(\|\rho u_t\|_{L^q}+ \|\rho u\cdot\na u\|_{L^q} + \||H||\na H|\|_{L^q} \right) \\
&\le C \left(\|\n u_t\|_{L^q} + \|\n u\|_{L^{2q}} \|\na u\|_{L^{2q}}+ \|H\|_{L^{2q}} \|\na H\|_{L^{2q}}
\right) \\
& \le C \|\n  u_t\|_{L^2}^{2(q-1)/(q^2-2)}\|\n u_t\|_{L^{q^2}}^{(q^2-2q)/(q^2-2)}\\
&\quad+ C\psi^\alpha\left(1+ \|\na^2 u\|_{L^{2}}^{1-1/q} +  \|\na^2 H\|_{L^{2}}^{1-1/q}  \right) \\
&\le  C  \left(\|\n^{1/2}  u_t\|_{L^2}^{2(q-1)/(q^2-2)}\|\na u_t\|_{L^{2}}^{(q^2-2q)/(q^2-2)}+\|\n^{1/2}  u_t\|_{L^2}\right)\\
&\quad+ C\psi^\alpha\left(1+  \|\na^2 u\|_{L^{2}}^{1-1/q} +  \|\na^2 H\|_{L^{2}}^{1-1/q}  \right),
\ea\ee
which together with  \eqref{gj3} and \eqref{gj6} implies that
\be \la{olv3.61} \ba \quad& \int_0^t\xl(\|\na^2u\|_{L^q}^{(q+1)/q}+\|\na P\|_{L^q}^{(q+1)/q}\xr)ds\\
\le  &C\int_0^t  s^{-\frac{q+1}{2q}} \left(s\|\n^{1/2} u_t\|_{L^2}^2\right)^{\frac{q^2-1}{ q(q^2-2)}} \left(s\|\na u_t\|^2_{L^2}\right)^{\frac{(q-2)(q+1)}{2 (q^2-2)}} ds \\
&+C\int_0^t   \|\n^{1/2}  u_t\|_{L^2}^{\frac{q+1}{q}}ds
+C\int_0^t \psi^\alpha \left( 1+ \|\na^2 u\|_{L^2}^{\frac{q^2-1}{q^2}}+ \|\na^2 H\|_{L^2}^{\frac{q^2-1}{q^2}}\right)ds \\
\le  &C\sup_{0\le s\le t} \left(s\|\n^{1/2} u_t\|_{L^2}^2\right)^{\frac{q^2-1}{ q(q^2-2)}}\int_0^t
s^{-\frac{q+1}{2q}}\left(s\|\na u_t\|^2_{L^2}\right)^{\frac{(q-2)(q+1)}{2 (q^2-2)}} ds \\
&+C\int_0^t \left( \psi^\alpha+\|\n^{1/2} u_t\|_{L^2}^2+ \|\na^2 u\|_{L^2}^{2}+ \|\na^2 H\|_{L^2}^2\right)ds \\
\le &C \exp\left\{C \int_0^t\psi^\alpha ds\right\}\left(1+\int_0^t \left( s^{-\frac{q^3+q^2-2q-2}{q^3+q^2-2q}}+s\|\na  u_t\|_{L^2}^2\right)ds\right)\\
\le& C \exp\left\{C \int_0^t\psi^\alpha ds\right\}
\ea\ee
and
\be\ba\la{olv3.62}
&\int_0^t \xl(s\|\na^2u\|_{L^q}^2+s\|\na P\|_{L^q}^2\xr)ds\\
\le & C \int_0^t  s \|\n^{1/2}  u_t\|_{L^2}^2ds +C\int_0^t  \left(s\|\n^{1/2}  u_t\|_{L^2}^2\right)^{2(q-1)/(q^2-2)}\left(s\|\na u_t\|_{L^{2}}^2\right)^{(q^2-2q)/(q^2-2)}ds\\
&+ C\int_0^t s\psi^\alpha\left( 1+ \|\na^2 u\|_{L^{2}}^{1-1/q} +  \|\na^2 H\|_{L^{2}}^{1-1/q}  \right)^2ds\\
\le & C \int_0^t  s \|\n^{1/2}  u_t\|_{L^2}^2ds +C\int_0^t s\|\na u_t\|_{L^{2}}^2 ds + C\int_0^t \left(\psi^\alpha+ s\|\na^2 u\|_{L^{2}}^2+  \|\na^2 H\|_{L^{2}}^2\right)ds\\
\le &C  \exp\left\{C \int_0^t\psi^\alpha ds\right\}.
\ea\ee
One thus obtains \eqref{gj7} from \eqref{olv3.61}--\eqref{olv3.62} and finishes the proof of Lemma \ref{le-3'}. \hfill $\Box$

\begin{lemma}\la{le-4}
Let $(\n,u,P,H)$ and $T_1$ be as in Lemma \ref{l3.01}. Then there exists a positive constant $\alpha>1$ such that for all $t\in (0, T_1]$,
\be\ba\la{gj8}
&\sup\limits_{0\le s\le t} \|\n \bar x^a \|_{ L^1\cap H^1\cap W^{1,q}}\le  \exp\left\{C\exp\left\{C\int_{0}^{t} \psi^{\alpha} ds\right\}\right\}.
\ea\ee
\end{lemma}

\pf First, it follows from Sobolev's inequality and \eqref{3.a2} that
for $0<\de<1,$
\be\ba\la{3.22}
\norm[L^{\infty}]{u\bar x^{-\de}}
&\le C(\de)\left(\norm[L^{4/\de}]{u\bar x^{-\de}}+ \norm[L^{3}]{\na (u\bar x^{-\de})}\right) \\ & \le C(\de)\left(\norm[L^{4/\de}]{u\bar x^{-\de}}+ \norm[L^{3}]{\na u}+\norm[L^{4/\de}]{u\bar x^{-\de}} \| \bar {x}^{-1}\na\bar{x} \|_{L^{12/(4-3\de)}} \right) \\
& \le C(\de) \left(\psi^\alpha+\norm[L^2]{\na^2 u} \right).
\ea\ee

One derives from \eqref{mhd}$_1$ that $\rho\bar{x}^a$ satisfies
\begin{equation}\label{A}
\partial_{t}(\rho\bar{x}^a)+ u \cdot\nabla(\rho\bar{x}^a)
-a\rho\bar{x}^{a} u \cdot\nabla\log\bar{x}=0,
\end{equation}
which along with \eqref{3.22} gives that
for any $r\in[2,q]$,
\be\ba \label{6.4}
\frac{d}{dt}\|\nabla(\rho\bar{x}^a)\|_{L^r} \leq &
C\left(1+\|\nabla u \|_{L^\infty}
+\|u\cdot\nabla\log\bar{x}\|_{L^\infty}\right)
\|\nabla(\rho\bar{x}^a)\|_{L^r} \\
& +C\|\rho\bar{x}^a\|_{L^\infty}\left(\||\nabla u ||\nabla\log\bar{x}|\|_{L^r}+\||u||\nabla^{2}\log\bar{x}|\|_{L^r}\right) \\ \leq & C\left(\psi^\alpha+\|\nabla^2 u \|_{L^2\cap L^q}\right)\|\nabla(\rho\bar{x}^a)\|_{L^r} \\
&+ C\|\rho\bar{x}^a\|_{L^\infty}\left(\|\nabla u \|_{L^r}
+\|u\bar{x}^{-\frac{2}{5}}\|_{L^{4r}}\|\bar{x}^{-\frac{3}{2}}\|_{L^{\frac{4r}{3}}}
\right) \\
\leq & C\left(\psi^\alpha+\|\nabla^2 u \|_{L^2\cap L^q}\right)
\left(1+\|\nabla(\rho\bar{x}^a)\|_{L^r}+\|\nabla(\rho\bar{x}^a)\|_{L^q}\right),
\ea\ee
where in   the last inequalities one has used  \eqref{igj1-2}.

Finally, using  \eqref{gj3}, \eqref{gj7}, \eqref{igj1}, \eqref{6.4}, and Gronwall's inequality, one thus gets \eqref{gj8} and completes the proof of Lemma \ref{le-4}.   \hfill $\Box$

Now, Proposition \ref{pro} is a direct consequence of  Lemmas \ref{l3.00}--\ref{le-4}.

\emph{Proof of Proposition \ref{pro}}.
It follows from \eqref{gj1}, \eqref{igj1}, \eqref{gj3}, and \eqref{gj8} that
\bnn\ba\psi(t)
&\le \exp\left\{C\exp\left\{C\int_{0}^{t} \psi^{\alpha} ds\right\}\right\}.
\ea\enn
Standard arguments yield  that for $M\triangleq e^{Ce}$ and $T_0\triangleq \min\{T_1,(CM^\alpha)^{-1}\}$,
\bnn \sup\limits_{0\le t\le T_0}\psi(t)\le M,\enn
which together with \eqref{gj3}, \eqref{gj6}, \eqref{gj10'}, and \eqref{gj7} gives \eqref{o1}. The proof of Proposition \ref{pro} is thus completed.   \hfill $\Box$

\section{Proof of Theorem \ref{t1}}

With the a priori estimates in Section 3 at hand, it is a position to prove Theorem 1.1.

{\it Proof of Theorem  \ref{t1}.}
Let $(\n_{0},u_{0},H_{0})$ be as in Theorem \ref{t1}. Without loss of generality, we assume that the initial density $\n_0$ satisfies
\bnn \int_{\rr} \n_0dx=1,\enn
which implies that there exists a positive constant $N_0$ such that  \be\la{oi3.8} \int_{B_{N_0}} \n_0 dx\ge \frac34\int_{\rr}\n_0dx=\frac34.\ee
We construct
$\n_{0}^{R}=\hat\n_{0}^{R}+R^{-1}e^{-|x|^2} $, where $0\le\hat\n_{0}^{R}\in  C^\infty_0(\rr)$ satisfies
\be\la{bci0}
\begin{cases}
\int_{B_{N_0}}\hat\n^R_0dx\ge 1/2,\\
\bar x^a \hat\n_{0}^{R}\rightarrow \bar x^a \n_{0}\quad {\rm in}\,\, L^1(\rr)\cap H^{1}(\rr)\cap W^{1,q}(\rr) ,\,\,{\rm as}\,\,R\rightarrow\infty. \end{cases}
\ee
Noting that $H_0\bar{x}^{a/2}\in L^2(\rr)$ and $\na H_0\in L^2(\rr)$, we choose $H_0^R\in \{w\in C^\infty_0(\O)~|~\div w=0\}$ satisfying
\be\ba\la{lv6.1}
H_0^R\bar x^{a/2} \rightarrow  H_0\bar x^{a/2} ,\quad \na H_{0}^{R}\rightarrow \na H_{0}\quad
{\rm in}\,\, L^2(\rr),\quad {\rm as}\,\,R\rightarrow\infty.
\ea\ee
Since $\na u_0\in L^2(\rr),$  we select $v^R_i\in C^\infty_0(\O)~(i=1,2 )$ such that for $i=1,2,$
\be\la{bci3}
\lim\limits_{R\rightarrow \infty}\|v^R_i-\pa_iu_0\|_{L^2(\rr)}=0.
\ee

We consider the unique smooth solution $u_0^R$ of the following elliptic problem:
\be\la{bbi2}
\begin{cases}
-\lap u_{0}^{R}+\n_0^R u_0^R+\na P^{R}_0 =\sqrt{\n_{0}^{R}} h^R-  \p_iv^R_i ,& {\rm in} \,\,  B_{R},\\
\div u_{0}^{R} =0,\,\, \,& {\rm in}~ B_{R},\\
u_{0}^{R} =0,\,\, \,& {\rm on}~\partial B_{R},
\end{cases}
\ee
where $h^R=(\sqrt{\n_0}u_0)*j_{1/R}$ with $j_\de$ being the standard mollifying kernel of width $\de.$

Extending $u_{0}^{R} $ to $\rr$ by defining $0$ outside $B_{R}$ and denoting it by $\tilde{u}_{0}^{R}$, we claim that
\be\la{3.7i4}
\lim\limits_{R\rightarrow \infty}\left(\|\na( \tilde{u}_{0}^{R}-u_0)\|_{L^2(\rr)}+\|\sqrt{\n_0^R}  \tilde{u}_{0}^{R}-\sqrt{\n_0}u_0 \|_{L^2(\rr)}\right)=0.
\ee
In fact, it is easy to find that $\tilde{u}_0^R$ is also a solution of \eqref{bbi2} in $\rr$. Multiplying \eqref{bbi2} by $\tilde{u}_0^R$ and integrating the resulting equation over $\rr$ lead to
\bnn\ba
&\int_{\rr} \n_0^R |\tilde{u}_0^R|^2dx +\int_{\rr} |\na \tilde{u}_0^R|^2dx\\
&\le \|\sqrt{\n_0^R} \tilde{u}_0^R\|_{L^2(\O)} \|h^R\|_{L^2(\O)} +C\|v_i^R\|_{L^2(\O)}\|\p_i \tilde{u}^R_0\|_{L^2(\O)}\\
&\le \frac{1}{2} \|\na \tilde{u}_0^R\|_{L^2(\O)}^2+ \frac{1}{2} \int_{\O} \n_0^R |\tilde{u}_0^R|^2dx+C\|h^R\|_{L^2(\O)}^2+C\|v_i^R\|_{L^2(\O)}^2,
\ea\enn
which implies \be \la{2.i9-4} \int_{\rr} \n_0^R|\tilde{u}_0^R|^2dx+\int_{\rr} |\na \tilde{u}_0^R|^2dx \le C \ee for some $C$ independent of $R.$
This together with \eqref{bci0} yields that there exist a subsequence $R_j\rightarrow \infty$ and a function $\tilde{u}_0\in\{\tilde{u}_0\in H^1_{\rm loc}(\rr)|\sqrt{\n_0}\tilde{u}_0\in L^2(\rr), \na \tilde{u}_0\in L^2(\rr)\}$ such that
\be\la{bci9}
\begin{cases}
\sqrt{\n^{R_j}_0}\tilde{u}^{R_j}_0 \rightharpoonup \sqrt{\n_0} \tilde{u}_0 \mbox{ weakly in } L^2(\rr) ,\\
\na \tilde{u}_0^{R_j}\rightharpoonup \na \tilde{u}_0 \mbox{ weakly in } L^2(\rr).
\end{cases}
\ee
Next, we will show
\be\la{bai1} \tilde{u}_0=u_0.\ee
Indeed, multiplying \eqref{bbi2} by a test function $\pi\in C_0^\infty(\rr)$ with $\div\pi=0$, it holds that
\be\la{ibai1}
\int_{\rr}\p_i(\tilde{u}_0^{R_j}-u_0)\cdot\p_i\pi dx+\int_{\rr} \sqrt{\n^{R_j}_0}(\sqrt{\n^{R_j}_0}\tilde{u}^{R_j}_0-h^{R_j})\cdot\pi dx=0.\ee
Let $R_j\rightarrow \infty$, it follows from \eqref{bci0}, \eqref{bci3}, and \eqref{bci9} that
\be\la{ibai2} \int_{\rr}\p_i(\ti u_0-u_0)\cdot\p_i\pi dx+\int_{\rr} \n_0(\ti u_0-u_0)\cdot\pi dx=0,\ee
which implies \eqref{bai1}.

Furthermore, multiplying  \eqref{bbi2} by $\tilde{u}_0^{R_j}$ and integrating the resulting equation over $\rr$, by the same   arguments as \eqref{ibai2}, we have
\bnn\ba
\lim \limits_{R_j\rightarrow \infty}\int_{\rr}\left( |\na \tilde{u}_0^{R_j}|^2  + \n_0^{R_j}|\tilde{u}_0^{R_j}|^2\right)dx = \int_{\rr}\left( |\na u_0 |^2  + \n_0 |u_0 |^2\right)dx,
\ea\enn
which combined with \eqref{bci9} leads to
\bnn
\lim\limits_{R_j\rightarrow \infty}\int_{\rr}|\na \tilde{u}_0^{R_j} |^2dx=\int_{\rr}|\na \tilde{u}_0 |^2dx,\,\,\lim\limits_{R_j\rightarrow \infty}\int_{\rr} \n_0^{R_j}   |\tilde{u}_0^{R_j} |^2dx=\int_{\rr}\n_0 | \tilde{u}_0 |^2dx.
\enn
This, along with  \eqref{bai1} and \eqref{bci9}, gives \eqref{3.7i4}.

Hence, by virtue of Lemma \ref{th0}, the initial-boundary-value problem  \eqref{mhd} and \eqref{ib2} with the initial data $(\n_0^R,u_0^R,H_0^R)$ has  a classical solution $(\n^{R},u^{R},P^{R},H^{R})$ on $B_{R}\times [0,T_R].$  Moreover, Proposition \ref{pro} shows that there exists a $T_0$ independent of $R$ such that \eqref{o1} holds for $(\n^{R},u^{R},P^{R},H^{R})$.

For simplicity, in what follows, we denote
\bnn L^p=L^p(\rr),\quad W^{k,p}=W^{k,p}(\rr).\enn

Extending $( \n^{R},u^{R},P^{R},H^{R})$ by zero on $\rr\setminus B_{R}$ and denoting it by
$$\xl(\nr\triangleq \vp_R\n^R,  \tilde{u}^R, \tilde{P}^{R}, \tilde{H}^R\xr)$$
with $\vp_R$ satisfying \eqref{vp1}. First, \eqref{o1} leads to
\be\ba\la{kq1}
&\sup\limits_{0\le t\le T_0}\left(\|\sqrt{\ti\n^R }  \tilde{u}^R\|_{L^2}+\|\na  \tilde{u}^R\|_{L^2}+\|\na\tilde{H}^R\|_{L^2}+\|\tilde{H}^R\bar x^{a/2} \|_{L^2}\right)\\
&\le \sup\limits_{0\le t\le T_0}\left(\|\sqrt{\n^R }  u^R\|_{L^2(\O)}+\|\na  u^R \|_{L^2(\O)}+\|\na  H^R \|_{L^2(\O)}+\|H^R \bar x^{a/2} \|_{L^2(\O)}\right)\\ &\le C ,
\ea\ee
and
\be\ba
&\sup\limits_{0\le t\le T_0}\|\ti\n^R\bar x^a\|_{L^1\cap L^\infty}\le C. \ea\ee
Similarly, it follows from \eqref{o1} that for $q>2$,
\be\ba
&\sup\limits_{0\le t\le T_0}t^{1/2}\left(\|\sqrt{\ti\n^R }  \tilde{u}^R_t\|_{L^2}+\|\na^2\tilde{u}^R\|_{L^2}+\|\na^2\tilde{H}^R \|_{L^2 }+\|\tilde{H}^R_t\|_{L^2 }\right)\\
&+\int_0^{T_0}\left(\|\sqrt{\ti\n^R}\ti u^R_t\|_{L^2}^2+\|\na^2\ti u^R\|_{L^2 }^2 +\|\ti H^R_t\|_{L^2 }^2+ \| \Delta \ti H^R\|_{L^2 }^2+\|\na \ti H^R\bar x^{a/2}\|_{L^2}^2 \right) dt\\
&+\int_0^{T_0} \left(\norm[L^{q}]{\nabla^2\ti u^R}^{(q+1)/q}+t\|\na^2\ti u^R\|_{L^q }^{2}+t\|\na \ti u^R_t\|_{L^2 }^2+t\|\na\ti H^R_t\|_{L^2 }^2 \right) dt\\
&\le C.
\ea\ee

Next, for $p\in[2,q]$, we obtain from \eqref{o1} and \eqref{gj8} that
\be\ba
\sup\limits_{0\le t\le T_0} \|\na (\ti\n^R\bar x^a) \|_{L^p }
&\le C\sup\limits_{0\le t\le T_0}\left( \|\na (\n^R\bar x^a)\|_{L^p(\O)}+R^{-1}\| \n^R\bar x^a \|_{L^p(\O)}\right)\\
&\le C \sup\limits_{0\le t\le T_0} \| \n^R\bar x^a\|_{H^1(\O)\cap W^{1,p}(\O)}\le C,
\ea\ee
which together with \eqref{3.22} and \eqref{o1} yields
\be\ba
\int_0^{T_0}\|\bar x\ti\n^R_t\|^2_{L^p }dt&\le C\int_0^{T_0} \|\bar x |u^R||\na \n^R| \|^2_{L^p(\O)} )dt\\ &\le C\int_0^{T_0} \|\bar x^{1-a} u^R\|_{L^\infty(\O)}^2\|\bar x^a\na \n^R  \|^2_{L^p(\O)}dt \\&\le C.
\ea\ee

By virtue of the same arguments as those of \eqref{gj10'} and \eqref{gj7}, one gets
\be\la{kq2}\ba  &\sup\limits_{0\le t\le T_0}t^{1/2}\|\na \ti P^R\|_{L^2 }+\int_0^{T_0}\left(\|\na \ti P^R\|_{L^2 }^2+\|\na \ti P^R\|_{L^q }^{(1+q)/q}\right)dt\le C.\ea\ee

With the estimates \eqref{kq1}--\eqref{kq2} at hand, we find that the sequence
$(\ti\n^R,\ti u^R,\ti P^R, \ti H^R)$ converges, up to the extraction of subsequences, to some limit $(\n,u,P, H)$ in the
obvious weak sense, that is, as $R\rightarrow \infty,$ we have
\be \la{kq3}\ti\n^R\bar x\rightarrow \n\bar  x, \mbox{ in } C(\overline{ B_N}\times [0,T_0]), \mbox{ for any }N>0,\ee
\be \ti\n^R\bar x^a\rightharpoonup  \n \bar x^a,\mbox{ weakly * in }L^\infty(0,T_0; H^1 \cap W^{1,q}),\ee
\be\ti H^R \bar x^{a/2}\rightharpoonup H\bar x^{a/2},\mbox{ weakly * in }L^\infty(0,T_0; L^2),\ee
\be  \sqrt{\ti\n^R} \ti u^R\rightharpoonup \sqrt{\n } u,
\,\, \na \ti u^R\rightharpoonup \na u ,\,\,
\na \ti H^R\rightharpoonup \na H,  \mbox{ weakly * in }L^\infty(0,T_0; L^2),\ee
\be \na^2 \ti u^R\rightharpoonup \na^2 u ,\,\,  \na  \ti P^R\rightharpoonup  \na P ,\mbox{ weakly  in }L^{\frac{q+1}{q}}(0,T_0; L^q)\cap L^2(\rr\times (0,T_0)),\ee
\be   \ti H_t^R\rightharpoonup  H_t ,\,\, \na  \ti H^R\bar x^{a/2}\rightharpoonup \na H \bar x^{a/2},\,\, \na^2 \ti H^R\rightharpoonup \na^2 H ,\mbox{ weakly  in } L^2(\rr\times (0,T_0)),\ee
\be   \sqrt{t }\na^2 \ti u^R\rightharpoonup  \sqrt{t }\na^2 u ,\mbox{ weakly  in }L^2(0,T_0; L^q),
\mbox{ weakly * in }L^\infty(0,T_0; L^2),\ee
\be \sqrt{t } \sqrt{\ti \n^R} \ti u^R_t\rightharpoonup  \sqrt{t }\sqrt{ \n } u_t,\,\,  \sqrt{t }\na  \ti P^R\rightharpoonup  \sqrt{t }\na P ,\mbox{ weakly * in }L^\infty(0,T_0; L^2),\ee
\be  \sqrt{t }\ti H^R_t\rightharpoonup  \sqrt{t }H_t, \,\,\sqrt{t }\Delta\ti H^R\rightharpoonup  \sqrt{t}\Delta H,\,\,\mbox{ weakly * in }L^\infty(0,T_0; L^2),
\ee
\be
\sqrt{ t} \na  \ti u^R_t\rightharpoonup \sqrt{t }\na  u_t,\quad \sqrt{t}\na  \ti H^R_t\rightharpoonup \sqrt{t }\na  H_t,\,\,\mbox{ weakly  in } L^2(\rr\times (0,T_0)),
\ee
with \be \la{kq4}  \n \bar x^a\in L^\infty(0,T_0; L^1), \quad \inf\limits_{0\le t\le T_0}\int_{B_{2N_0}}\n(x,t)dx\ge \frac14. \ee

Then letting $R\rightarrow \infty$, standard arguments together with  \eqref{kq3}--\eqref{kq4} show that
$(\n,u,P,H)$ is a strong solution of \eqref{mhd}-\eqref{n4} on $\rr\times (0,T_0]$ satisfying \eqref{1.10} and \eqref{l1.2}. Indeed, the existence of a pressure $P$ follows immediately from the \eqref{mhd}$_2$ and \eqref{mhd}$_4$  by a classical consideration. The  proof of the existence part of Theorem \ref{t1} is finished.

It remains only  to prove the uniqueness of the strong solutions  satisfying \eqref{1.10} and \eqref{l1.2}.
Let $(\rho,u,P,H)$ and $(\bar\rho,\bar u,\bar P, \bar H)$ be two strong solutions satisfying \eqref{1.10} and \eqref{l1.2} with the same initial data, and denote
$$\Theta\triangleq\rho-\bar\rho,~U\triangleq u-\bar u,~\Phi\triangleq H-\bar H.$$

First, subtracting the mass equation satisfied by $(\rho,u,P,H)$ and $(\bar\rho,\bar u,\bar P,\bar H)$ gives
\be\la{5.2}
\Theta_t+\bar u\cdot\nabla\Theta+U\cdot\nabla \rho= 0.
\ee
Multiplying \eqref{5.2} by $2\Theta\bar{x}^{2r}$ for $r\in (1,\tilde a)$ with $\tilde a=\min\{2,a \}$, and integrating by parts yield
\bnn\ba
&\frac{d}{dt}\int|\Theta\bar{x}^{r}|^{2}dx\\
&\le C \norm[L^{\infty}]{ \bar{u} \bar{x}^{-1/2}} \norm[L^{2}]{\Theta\bar{x}^{r}}^{2} +C\norm[L^{2}]{\Theta\bar{x}^{r}}\norm[L^{2q/((q-2)(\tilde a-r))}]{ U  \bar{x}^{-(\tilde a-r)}}\norm[L^{2q/(q-(q-2)(\tilde a-r))}]
{\bar{x}^{\tilde a} \na \n} \\
&\le C\left( 1+\norm[ W^{1,q}]{\na\bar u}  \right)\norm[L^{2}]{\Theta\bar{x}^{r}}^{2}+C\norm[L^{2}]{\Theta\bar{x}^r}
\left(\norm[L^{2}]{\nabla U}+\norm[L^{2}]{\sqrt{\n} U}\right)
\ea\enn
due to Sobolev's inequality, \eqref{l1.2}, \eqref{3.v2}, and \eqref{3.22}. This combined with Gronwall's inequality shows that  for all $0\le t\le T_{0}$,
\be\ba\la{5.1}
\norm[L^{2}]{\Theta\bar{x}^{r}}
\le & C\int_{0}^{t}\left(\norm[L^{2}]{\nabla U}+\norm[L^{2}]{\sqrt{\n} U}\right)ds.
\ea\ee

Next, subtracting the momentum and magnetic equations satisfied by $(\rho,u,P,H)$ and $(\bar\rho,\bar u,\bar P,\bar H)$ leads to
\be\ba\la{5.5}
\rho U_t + \rho u\cdot\nabla U -\mu\lap U
= & - \rho U\cdot\nabla\bar u-\Theta(\bar u_t+\bar u\cdot\nabla\bar u)-\na ( P -\bar P) \\& -\frac12\na\left(|H|^2-|\bar H|^2\right)+H\cdot\na \Phi+\Phi\cdot\na \bar H,
\ea\ee
and
\be\la{hlv7.1}
\Phi_t-\nu\Delta \Phi=H\cdot\na U+\Phi\cdot\na\bar u-u\cdot\na\Phi-U\cdot\na\bar H.
\ee
Multiplying (\ref{5.5}) and \eqref{hlv7.1} by $U$ and $\Phi$, respectively, and adding the resulting equations together, we obtain after integration by parts that
\be\ba\la{5.6}
& \frac{d}{dt}\int \left(\rho |U|^2+|\Phi|^2\right)dx  +  \int \left(\mu|\na U|^2+\nu|\na \Phi|^2\right) dx\\
& \le   C \norm[L^{\infty}]{\nabla \bar{u}} \int \left(\n |U|^2+|\Phi|^2\right)dx+C\int|\Theta|| U| \left(| \bar u_{t}|+ | \bar u ||\na\bar{u}|\right) dx  \\
&\quad -\int\Phi\cdot\na U\cdot\bar Hdx -\int U\cdot\na \bar H\cdot\Phi dx\\
& \triangleq   C  \norm[L^{\infty}]{\nabla \bar{u}} \int \left(\n |U|^2+|\Phi|^2\right)dx+\sum_{i=1}^3K_i.
\ea\ee

We first estimate $K_1$.  H\"older's inequality combined with \eqref{l1.2}, \eqref{3.i2}, \eqref{o1}, and \eqref{5.1} yields that for $r\in (1,\tilde a),$
\be\ba\la{hlv7.4}
K_1& \le C\norm[L^2]{\Theta \bar{x}^{r}}\norm[L^{4}] {U \bar{x}^{-r/2}}\left( \norm[L^{4}]{\bar u_{t}\bar{x}^{-r/2}}+\norm[L^{\infty}]{\na \bar{u}} \norm[L^{4}] {\bar u\bar{x}^{-r/2}}\right) \\
&\le C(\ve)\left( \norm[L^{2}] {\sqrt{\bar{\n}}\bar{u}_{t}}^{2}+ \norm[L^{2}] {\na\bar{u}_{t}}^{2}+\norm[L^{\infty}]{\na \bar{u}}^{2}\right)\norm[L^2]{\Theta \bar{x}^{r}}^2\\ & \quad+\ve\left(\norm[L^{2}]{\sqrt{\rho} U}^{2}+ \norm[L^{2}]{\na U}^2\right)\\
&\le C(\ve)\left(1+t\norm[L^{2}] {\nabla  {\bar{u}_t}}^{2}+t\norm[L^{q}]{ \na^{2} \bar{u}}^{2} \right)\int_{0}^{t} \left(\norm[L^{2}]{\nabla U}^{2}+\norm[L^{2}]{\sqrt{\n} U}^{2}\right)ds\\
& \quad+\ve \left(\norm[L^{2}]{\sqrt{\rho} U}^{2}+ \norm[L^{2}]{\na U}^2\right).
\ea\ee
For the term $K_2$, we derive from Gagliardo-Nirenberg inequality and \eqref{gj2} that
\be\ba\la{hlv7.6}
K_2 \le&C \|\bar H\|_{L^4} \|\Phi\|_{L^4}\|\na U\|_{L^2}\\
\le &\ep\|\na U\|_{L^2}^2+\ep\|\na \Phi\|_{L^2}^2+C(\ep)\|\Phi\|_{L^2}^2.
\ea\ee
The last term $K_3$ can be estimated as follows
\be\ba\la{hlv7.7}
K_3\le &C\|U\bar x^{-a}\|_{L^4}\||\na \bar H|^{1/2}\bar x^{a}\|_{L^4}\||\na \bar H|^{1/2}\|_{L^4}\|\Phi\|_{L^4}\\
\le & C\left(\|\sqrt{\n}U\|_{L^2}+\|\na U\|_{L^2}\right)\|\na \bar H\bar x^{a/2}\|_{L^2}^{1/2}\|\Phi\|_{L^4}\\
\le & \ve\left(\|\sqrt{\n}U\|_{L^2}^2+\|\na U\|_{L^2}^2\right)+C(\ve)\|\na \bar H\bar x^{a/2}\|_{L^2}\|\Phi\|_{L^4}^2\\
\le & \ve\left(\|\sqrt{\n}U\|_{L^2}^2+\|\na U\|_{L^2}^2\right)+\ve\|\na \Phi\|_{L^2}^2+C(\ve)\|\na \bar H\bar x^{a/2}\|_{L^2}^2\|\Phi\|_{L^2}^2
\ea\ee
owing to \eqref{l1.2}, \eqref{3.i2}, and \eqref{o1}.

Denoting
$$
G(t)\triangleq \norm[L^{2}]{\sqrt{\rho} U}^{2}+\|\Phi\|_{L^2}^2+\int_0^t\left( \|\na U\|_{L^2}^2+\|\na \Phi\|_{L^2}^2+\|\sqrt{\n}U\|_{L^2}^2\right)ds,
$$
then substituting \eqref{hlv7.4}--\eqref{hlv7.7} into \eqref{5.6} and choosing $\ep$ suitably small lead  to
$$
G'(t)\le C\left(1 +\|\na \bar u\|_{L^\infty}+\|\na \bar H \bar x^{a/2}\|_{L^2}^2+t\|\na \bar u_t\|_{L^2}^2+t\|\na^2 u\|_{L^q}^2\right)G(t),
$$
which together with Gronwall's inequality and \eqref{1.10} implies $G(t)=0$. Hence, $U(x,t)=0$ and $\Phi(x,t)=0$ for almost everywhere $(x, t)\in\mr^2\times(0, T)$. Finally, one can deduce from \eqref{5.1} that $\Theta=0$ for almost everywhere $(x, t)\in\mr^2\times(0, T)$.
The proof of Theorem 1.1 is completed.   \hfill $\Box$

\end{document}